\input amstex\documentstyle{amsppt}  
\pagewidth{12.5cm}\pageheight{19cm}\magnification\magstep1  
\topmatter
\title On induction of class functions\endtitle
\author G. Lusztig\endauthor
\address{Department of Mathematics, M.I.T., Cambridge, MA 02139}\endaddress
\thanks{Supported by NSF grant DMS-1855773.}\endthanks
\endtopmatter   
\document

\define\si{\sim}

\define\sqc{\sqcup}

\define\bY{\bar Y}

\define\op{\oplus}
   
\define\part{\partial}
\define\emp{\emptyset}

\define\m{\mapsto}

\define\sub{\subset}    
\define\bxt{\boxtimes}
\define\T{\times}
\define\ti{\tilde}
\define\nl{\newline}
\redefine\i{^{-1}}

\define\un{\underline}

\define\bbq{\bar{\QQ}_l}

\define\tr{\text{\rm tr}}

\redefine\b{\beta}
\redefine\c{\chi}
\define\g{\gamma}

\define\et{\eta}

\define\p{\pi}
\define\ph{\phi}
\define\ps{\psi}

\define\s{\sigma}
\redefine\t{\tau}
\define\th{\theta}

\define\z{\zeta}
\define\x{\xi}

\redefine\G{\Gamma}

\define\Si{\Sigma}
\define\Th{\Theta}
\redefine\L{\Lambda}

\define\kk{\bold k}

\define\QQ{\bold Q}

\define\ca{\Cal A}

\define\ce{\Cal E}
\define\cf{\Cal F}

\define\ck{\Cal K}
\define\cl{\Cal L}

\define\car{\Cal R}

\define\cz{\Cal Z}
\define\cx{\Cal X}
\define\cy{\Cal Y}

\define\tg{\ti g}

\define\tm{\ti m}

\define\tu{\ti u}
\define\tv{\ti v}

\define\tG{\ti G}

\define\tL{\ti L}
\define\tM{\ti M}

\define\tP{\ti P}

\define\tS{\ti S}

\define\tX{\ti X}
\define\tY{\ti Y}
\define\ttY{\ti{\ti Y}}

\define\sha{\sharp}

\define\tce{\ti\ce}

\head Introduction \endhead
\subhead 0.1\endsubhead
Let G be a connected reductive algebraic group over $\kk$, an algebraic closure of the finite field $F_q$
with $q$ elements, with a fixed $F_q$-rational structure whose Frobenius map is denoted by $F:G@>>>G$.

Let $\L(G)$ be the set of subgroups $M$ of $G$ such that $M$ is a Levi subgroup of a parabolic subgroup of $G$; 
for $M\in\L(G)$ let $\Pi(M)$ be the set of parabolic subgroups $P$ of $G$ for which $M$ is a Levi subgroup.
Assume that $M\in\L(G)$ is defined over $F_q$ and that $P\in\Pi(M)$ (not necessarily defined over $F_q$).
Let $\ck(G^F)$ (resp. $\ck(M^F)$) be the Grothendieck group of
representations of $G^F$ (resp. $M^F$) over an algebraic closure $\bbq$ of the $l$-adic numbers where $l$ is
a prime number not dividing $q$. (When $F$ acts on a set $X$ we denote by $X^F$ the fixed point set of 
$F:X@>>>X$.) Let $R_{M,P}^G:\ck(M^F)@>>>\ck(G^F)$ be the ``induction'' homomorphism defined in \cite{DL} (in the 
case where $M$ is a maximal torus) and in \cite{L76} (in the general case). 
Let $cl(G^F)$ (resp. $cl(M^F)$) be the $\bbq$-vector space of class functions $G^F@>>>\bbq$ 
(resp. $M^F@>>>\bbq$). By passage to characters, $R_{M,P}^G$ can be regarded as a $\bbq$-linear map
$R_{M,P}^G:cl(M^F)@>>>cl(G^F)$. In \cite{L76} it was conjectured that 

(a) {\it $R_{M,P}^G$ is independent of the choice of $P$.}
\nl
(At that time it was already known from \cite{DL} that (a) holds when $M$ is a maximal torus of $G$, so that in
that case, the notation $R_M^G$ can be used instead of $R_{M,P}^G$.) As noted in
\cite{L76}, Deligne had an argument to prove (a) for any $M$ provided that $q\gg0$ (but
his proof has not been published). In \cite{L90, 8.13} a proof of (a) for $q\gg0$ was given
which was based on the theory of character sheaves and thus was quite different from Deligne's proof. (In 
{\it loc. cit.} there is also an assumption on the characteristic
$p$ of $\kk$, but that assumption can be removed in view of the cleanness result for character sheaves in 
\cite{L12}.) In \cite{BM} it is proved that (a) holds assuming only that $q>2$.

In this paper we define a $\bbq$-linear map $\car_M^G:cl(M^F)@>>>cl(G^F)$ with no restriction on $q$ (see 1.7) 
and we show that 

(b) {\it if $q\gg0$ we have $\car_M^G=R_{M,P}^G$ for any $P\in\Pi(M)$.}
\nl
(See 1.9 and \S2).
We expect that the results of \cite{L90} quoted in this paper hold without restriction on $q$ and, as a
consequence, that (b) holds without restriction on $q$.

The definition of $\car_M^G$ is in terms of intersection cohomology; it relies on ideas of \cite{L84}. The
proof of (b) relies on the results of \cite{L90} connecting representations of $G^F$ with the character sheaves
on $G$. In \S3 we show (based on results of \cite{L90}) that if $D$ is an $F$-stable conjugacy class of $G^F$
then the function on $G^F$ which is $1$ on $D^F$ and $0$ on $G^F-D^F$ is a linear combination of characters of
$R_T^G(\th)$ for various $F$-stable maximal tori of $G$ and various characters $\th$ of $T^F$. (This was
conjectured in \cite{L78}.)

\subhead 0.2\endsubhead
{\it Notation.} Let $\nu_G$ be the dimension of the flag manifold of $G$. Let $\cz_G$ be the centre of $G$.
For $M\in\L(G)$ let $N_GM$ be the normalizer of $M$ in $G$. 
For $g\in G$ we have $g=g_sg_u=g_ug_s$ where $g_s\in G$ is semisimple and $g_u\in G$ is unipotent.
For $s\in G$ semisimple we write $G_s$ for the centralizer of $s$ in $G$. For $g\in G$ let $H_G(g)$ be the 
smallest subgroup in $\L(G)$ that contains $G_{g_s}^0$. If $G'$ is a subgroup of $G$, let $Z_G(G')$ be the 
centralizer of $G'$ in $G$. Let $G_{der}$ be the derived subgroup of $G$.

Let $X$ be an algebraic variety over $\kk$. Let $ls(X)$ be the collection of $\bbq$-local systems on $X$. If $H$ 
is a connected algebraic group acting on $X$ we denote by $ls_H(X)$ the collection of $H$-equivariant
$\bbq$-local systems on $X$. Let $Y$ be a locally closed, smooth, irreducible subvariety of $X$ and let 
$\ce\in ls(Y)$. Then
$\ce$ extends canonically as an intersection cohomology complex to the closure $\bY$ of $Y$ and to $X$ by $0$ on 
$X-\bY$; the resulting complex on $X$ is denoted by $\ce^{\sha}$. Assume now that $X$ is defined over $F_q$ with 
Frobenius map $F:X@>>>X$ and $Y$ above is $F$-stable. Assume that $F^*\ce\cong\ce$ and we are given an isomorphism
$\ph:F^*\ce@>\si>>\ce$. This induces an isomorphism $\ph^\sha: F^*\ce^\sha@>\si>>\ce^\sha$. 
Let $\c_{\ce,\ph}:X^F@>>>\bbq$ be the function whose value at $y\in Y^F$ is the trace of $\ph$ on the stalk of
$\ce$ at $y$ and which is zero on $X^F-Y^F$ 
Let $\c_{\ce^\sha,\ph^\sha}:X^F@>>>\bbq$ be the function whose value at 
$x\in X^F$ is the alternating sum of traces of the linear maps induced by $\ph^\sha$ on the stalks at $x$ of the 
cohomology sheaves of $\ce^\sha$. We have $\c_{\ce,\ph}|Y^F=\c_{\ce^\sha,\ph^\sha}|Y^F$.

\head 1. The definition of $\car_M^G$\endhead
\subhead 1.1\endsubhead
A subset $S$ of $G$ is said to be an isolated stratum if $S$ is the inverse image of an isolated conjugacy
class (see \cite{L84, 2.6}) of $G/\cz^0_G$ under the obvious map $G@>>>G/\cz^0_G$.
Let $A_G$ be the set of pairs $(L,S)$ where $L\in\L(G)$ and $S$ is an isolated stratum of $L$.
For $(L,S)\in A_G$ we set $S_r^G=\{g\in S;H_G(g)=L\}=\{g\in S;G_{g_s}^0\sub L\}$ and
$Y^G_{L,S}=\cup_{x\in G}xS_r^Gx\i$. Then 

(a) {\it $Y^G_{L,S}$ is a smooth locally closed irreducible subvariety of $G$ of dimension 
$2\nu_G-2\nu_L+\dim S$.}
\nl
(see \cite{L84, 3.1}). The subsets $Y^G_{L,S}$ are called the strata of $G$. We have 
$$G=\sqc_{(L,S)\in A_G\text{ up to G-conjugacy}}Y^G_{L,S}.$$
 Note that an isolated stratum $S$ of $G$ is the 
stratum $Y^G_{G,S}$.

\subhead 1.2\endsubhead
Let $M\in\L(G)$. Let $Y'$ be a stratum of $M$. We shall associate to $Y'$ a stratum $Y$ of $G$ as follows.
 We set $Y'_r=\cup_{x\in M}xS_r^Gx\i$, $Y=Y^G_{L,S}$ where 
$(L,S)\in A_M$ is such that $Y'=Y^M_{L,S}=\cup_{x\in M}xS_r^Mx\i$. (We have also $(L,S)\in A_G$.) Now $Y'_r$ and 
$Y$ are independent of the choice of $(L,S)$. (Indeed, it is enough to show that if $m\in M$, then 
$(mSm\i)^G_r=mS^G_rm\i$; we use that $H_G(mgm\i)=mH_G(g)m\i$.) We have $Y'_r\ne\emp$, $Y'_r\sub Y$ and $Y$
is the unique stratum of $G$ that contains $Y'_r$. We have also $Y'_r\sub Y'$. (We use that $S_r^G\sub S_r^M$; 
indeed, if $g\in S$ and $H_G(g)=L$, then $G_{g_s}^0\sub L$ hence $M_{g_s}^0\sub L$ and $H_M(g)=L$.) We show that

(a) {\it $Y'_r$ is open in $Y'$.}
\nl
Let $(L,S)\in A_M$ be such that $Y'=Y^M_{L,S}$. Now $Y'$ is a locally trivial fibration over the variety of
all $M$-conjugates of $L$, via $g\m H_M(g)$; let $\b$ be the fibre of this map over $L$. It is enough to show 
that $Y'_r\cap\b$ is open in $\b$. We have $\b=\cup_{n\in N_ML/L}nS^M_rn\i$, 
$Y'_r\cap\b=\cup_{n\in N_ML/L}nS^G_rn\i$. It is enough to observe that $S^G_r$ is open in $S^M_r$; in fact it is 
open in $S$.

Let $\ttY=\{(g,x)\in G\T G;x\i gx\in Y'_r\}$. Define $\s:\ttY@>>>Y'_r$ by $\s(g,x)=x\i gx$. Let 
$\tY=\ttY/M=\{(g,xM)\in Y\T G/M;x\i gx\in Y'_r\}$. We show:

(b) {\it $\tY$ is a smooth, irreducible variety of dimension equal to $\dim Y$.}
\nl
Since $Y'_r$ is smooth, irreducible of dimension equal to $\dim Y'$ (see (a), 1.1(a)), we see that $\tY$ is 
smooth, irreducible of dimension $\dim G/M+\dim Y'=\dim Y$. This proves (b).

If $(g,xM)\in\tY$, we have $g\in Y$ (since $Y'_r\sub Y$). Define $\p:\tY@>>>Y$ by $(g,xM)\m g$. We show:

(c) {\it $\p$ is a finite unramified cover of $Y$ with fibres isomorphic to $c_G/c_M$ where 
$c_G=\{n\in N_GL,n\i Sn=S\}$, $c_M=\{n\in N_ML,n\i Sn=S\}$ and $(L,S)\in A_M$, $Y'=Y^M_{L,S}$.}
\nl
It is enough to show that if $g\in S^G_r$, then $\p\i(g)\cong c_G$. We can identify $\p\i(g)$ with
$\{xM\in G/M;x\i gx\in Y'_r\}$ hence also with $\{x\in G;\x\i gx\in S^G_r\}/\{x\in M;\x\i gx\in S^G_r\}$. It is 
enough to show that $\{x\in G;\x\i gx\in S^G_r\}=c_G$ (this would imply $\{x\in M;\x\i gx\in S^G_r\}=c_M$).
Let $x\in G$ be such that $x\i gx\in S^G_r$; then $L=H_G(x\i gx)=x\i H_G(g)x=x\i Lx$ and $x\in N_GL$. Let $S',g'$
be the image of $S,g$ in $L/\cz^0_L$. Now $Ad(x\i)$ induces an automorphism of $L/\cz^0_L$ which carries 
$S'$ to an isolated conjugacy class $S''$ and $g'\in S'$ to an element of $S'$; it follows that $S'=S''$ and
$Ad(x\i)S=S$, so that $x\in c_G$. Conversely, let $x\in N_GL$ be such that $x\i Sx=S$. Then $x\i gx\in S$, 
$H_G(x\i gx)=x\i H_G(g)x=x\i Lx=L$ and $x\i gx\in S^G_r$. This proves (c).

\subhead 1.3\endsubhead
We preserve the setup of 1.2. Let $\ce\in ls_M(Y')$. We define $j_{Y'}^Y(\ce)\in ls_G(Y)$. Note that 
$\s^*(\ce|Y'_r)\in ls_{G\T M}(\ttY)$ for the $G\T M$ action $(g_0,m):(g,x)\m(g_0gg_0\i,g_0xm\i)$ on $\ttY$. 
Hence $\s^*(\ce|Y'_r)=\s_1^*\ce_1$ where $\s_1:\ttY@>>>\tY$ is the obvious map and $\ce_1\in ls(\tY)$ is
well defined. We define $j_{Y'}^Y\ce=\p_*(\ce_1)\in ls_G(Y)$; this has rank equal to $c_G/c_M$ times the 
rank of $\ce$.

\subhead 1.4\endsubhead
We preserve the setup of 1.3. Let $(L,S)\in A_M$ be such that $Y'=Y^M_{L,S}$ and let $\ce_0\in ls_L(S)$. Then 
$j_S^{Y'}(\ce_0)\in ls(Y')$ and $j_S^Y(\ce_0)\in ls(Y)$ are defined as in 1.3. From the definition we see that

(a) $j_{Y'}^Y(j_S^{Y'}\ce_0)=j_S^Y(\ce_0)$.

\subhead 1.5\endsubhead
Let $\cf\in ls_G(Y_1)$ where $Y_1=Y^G_{L,S}$. We say that $\cf$ is {\it admissible} if it is 
irreducible and if $\cf$ is a direct summand of $j_S^{Y_1}(\cf_0)\in ls_G(Y_1)$ for some $\cf_0\in ls_L(S)$
which is cuspidal irreducible (see \cite{L84, 2.4}).  (This condition is independent
of the choice of $(L,S)$.) Let $\ca^G(Y_1)$ be the class of $G$-equivariant admissible local systems on $Y_1$.
We say that $Y_1$ is an admissible stratum if $\ca^G(Y_1)\ne\emp$.
In the setup of 1.3 we see (using 1.4(a)) that if $\cf\in\ca^M(Y')$ then $j_{Y'}^Y(\cf)$ is a (nonzero) direct sum
of objects of $\ca^G(Y)$. In particular, if $Y'$ is admissible (for $M$) then $Y$ is admissible (for $G$).

\subhead 1.6\endsubhead
We now assume that $M$ is defined over $F_q$. If $Y'$ in 1.2 is $F$-stable then $Y'_r$ in 1.2 is $F$-stable. 
Indeed, from $Y'=F(Y')$ and $Y'=Y^M_{L,S}$ we deduce $Y'=F(Y')=Y^M_{F(L),F(S)}$ hence $F(L)=mLm\i$, $F(S)=mSm\i$
for some $m\in M$. Replacing $L,S$ by an $M$-conjugate we can assume that $F(L)=L,F(S)=S$, so that
$$F(Y'_r)=\cup_{x\in M}xF(S)_r^Gx\i=\cup_{x\in M}xS^G_rx\i=Y'_r.$$
A similar argument shows that $Y$ in 1.2 is $F$-stable (alternatively, this holds since $Y$ is the unique
stratum of $G$ containing $Y'_r$ which is $F$-stable). Moreover, $\ttY,\tY,\p$ in 1.2 and $\s,\s_1$ in 1.3 are 
defined over $F_q$.
If now $\ce$ in 1.3 is such that $F^*\ce\cong\ce$, then 
$\s^*(\ce|Y'_r)$ and $\ce_1$ in 1.3 are isomorphic to their inverse image under $F$ hence
$\tce:=j_{Y'}^Y(\ce)$ satisfies $F^*\tce\cong\tce$; moreover any isomorphism $\ph:F^*\ce@>\si>>\ce$ (of local 
systems on $Y'$) induces an isomorphism $F^*(\s^*(\ce|Y'_r))@>\si>>\s^*(\ce|Y'_r)$ (of local systems on $\ttY$),
an isomorphism $F^*\ce_1@>\si>>\ce_1$ (of local systems on $\tY$) and an isomorphism
$\ti\ph:F^*\tce\cong\tce$ (of local systems on $Y$). Now $\ph,\ti\ph$ extend to isomorphisms 
$\ph^\sha:F^*\ce^\sha@>\si>>\ce^\sha$ (of complexes on $M$) and $\ti\ph^\sha:F^*\tce^\sha@>\si>>\tce^\sha$ (of 
complexes on $G$). Hence $\c_{\ce^\sha,\ph^\sha}:M^F@>>>\bbq$, $\c_{\tce^\sha,\ti\ph^\sha}:G^F@>>>\bbq$ are well 
defined. They are class functions on $M^F$ (resp. $G^F$).

\subhead 1.7\endsubhead
Let $Z_M$ be the set of all pairs $(Y',\ce)$ where 
$Y'$ is an $F$-stable admissible stratum of $M$ and $\ce\in ls_M(Y')$ is admissible (up to isomorphism) such 
that $F^*\ce\cong\ce$. 
For any $(Y',\ce)\in Z_M$ we denote by $\cl_{Y',\ce}$ the subspace of $cl(M^F)$ consisting of the class
functions $\c_{\ce^\sha,\ph^\sha}:M^F@>>>\bbq$ where $\ph:F^*\ce@>\si>>\ce$ is an isomorphism. A different 
choice for $\ph$ must be of the form $a\ph$ for some $a\in\bbq^*$ and we have
$\c_{\ce^\sha,(a\ph)^\sha}=a\c_{\ce^\sha,\ph^\sha}$ hence $\cl_{Y',\ce}$ is well defined line in $cl(M^F)$, 
independent of any choice. We have

(a) $cl(M^F)=\op_{(Y',\ce)\in Z_M}\cl_{Y',\ce}$.
\nl
A proof is given in \cite{L04, 26.5}. (Alternatively, instead of \cite{L04}, one can use \cite{L86, 25.2} 
complemented by \cite{L12}.) Let 

(b) $\car_M^G:cl(M^F)@>>>cl(G^F)$
\nl
be the linear map such that for any $(Y',\ce)\in Z_M$, the restriction of $\car_M^G$ to the line $\cl_{Y',\ce}$
sends $\c_{\ce^\sha,\ph^\sha}:M^F@>>>\bbq$ (where $\ph:F^*\ce@>\si>>\ce$ is an isomorphism) to 
$\c_{\tce^\sha,\ti\ph^\sha}:G^F@>>>\bbq$ (see 1.6). If $\ph$ is replaced by $a\ph$ with $a\in\bbq^*$, then 
$\ti\ph^\sha$ is replaced by $a\ti\ph^\sha$. Thus the linear map $\car_M^G$ is well defined (independent of 
choices).

\subhead 1.8\endsubhead
In the setup of 1.4 assume that $\ce_0\in ls_L(S)$ is cuspidal irreducible, that $L,S$
are defined over $F_q$ and that we are given an isomorphism $\ph_0:F^*\ce_0@>\si>>\ce_0$ of local systems on $S$.
Then $\c_{\ce_0^\sha,\ph_0^\sha}:L^F@>>>\bbq$ is well defined. Using 1.4(a) and the definitions we see that

(a) $\car_M^G(\car_L^M(\c_{\ce_0^\sha,\ph_0^\sha}))=\car_L^G(\c_{\ce_0^\sha,\ph_0^\sha}))$.

\subhead 1.9\endsubhead 
In the rest of this paper, unless otherwise specified, we assume that $q\gg0$, so that the results of \cite{L90} can be
applied. (As mentioned in 0.1 the assumption of {\it loc. cit.} on the characteristic of $\kk$, can now be 
removed.) We shall write $R_M^G:cl(M^F)@>>>cl(G^F)$ instead of $R_{M,P}^G$ (with $P\in\Pi(M)$). 

Let $(L,S)\in A_M,\ce_0,\ph_0$ be as in 1.8. The following result can be deduced from \cite{L90, 9.2}:

(a) $\car_L^M(\c_{\ce_0^\sha,\ph_0^\sha})=R_L^M(\c_{\ce_0^\sha,\ph_0^\sha})$.
\nl
Let $cl'(M^F)$ be the $\bbq$-subspace of $cl(M^F)$ generated by the elements 
$R_L^M(\c_{\ce_0^\sha,\ph_0^\sha})\in cl(M^F)$ for various $L,S,\ce_0,\ph_0$ as in (a). The following result 
will be proved in \S2.

(b) {\it We have $cl'(M^F)=cl(M^F)$.}
\nl
We can now prove 0.1(b). By (b), it is enough to show that if $L,S,\ce_0,\ph_0$ are as in (a) then

(c) $R_M^G(R_L^M(\c_{\ce_0^\sha,\ph_0^\sha}))=\car_M^G(R_L^M(\c_{\ce_0^\sha,\ph_0^\sha}))$.
\nl
By (a), the right hand side of (c) is equal to $\car_M^G(\car_L^M(\c_{\ce_0^\sha,\ph_0^\sha}))$ hence, by 1.8(a),
it is equal to $\car_L^G(\c_{\ce_0^\sha,\ph_0^\sha})$. We have $R_M^GR_{L}^M=R_L^G:cl(L^F)@>>>cl(G^F)$. (This 
is proved in \cite{L76, Cor.5} assuming that $L$ is a maximal torus of $G$; but the same proof works in general.)
Thus (c) is equivalent to the equality

$R_L^G(\c_{\ce_0^\sha,\ph_0^\sha})=\car_L^G(\c_{\ce_0^\sha,\ph_0^\sha})$
\nl
and this follows from (a) (with $M$ replaced by $G$). This proves 0.1(b).

\head 2. Proof of 1.9(b)\endhead
\subhead 2.1\endsubhead
Let $\tG$ be a reductive connected group over $\kk$ with an $F_q$-rational structure (with Frobenius map 
$F:\tG@>>>\tG$) such that $\tG_{der}$ is simply connected; assume that we are given a surjective 
homomorphism of algebraic groups $\t:\tG@>>>G$ defined over $F_q$ whose kernel $K$ is a central torus in $\tG$.
Then $\tM=\t\i(M)\in\L(\tG)$ is defined over $F_q$. Let $P\in\Pi(M)$ and let $V$ be the unipotent radical of $P$.
Then $\tP=\t\i(P)\in\Pi(\tM)$. Let $\ti V$ be the unipotent radical of $P$.
Let $X=\{gV\in G;g\i F(g)\in F(V)\}$, $\tX=\{\tg\ti V\in\tG;\tg\i F(\tg)\in F(\ti V)\}$.
Now $G^F\T M^F$ acts on $X$ by $(g_0,m_0):g\m g_0gm_0\i$ and this induces an action of $G^F\T M^F$
on the $l$-adic cohomology $H^i_c(X,\bbq)$. Similarly, $\tG^F\T\tM^F$ acts on $H^i_c(\tX,\bbq)$.
For $u\in G^F$ unipotent and for $u'\in M^F$ unipotent we
set $\g_{M,V}^G(u,u')=\sum_i(-1)^i\tr((u,u'),H^i_c(X,\bbq))$. Similarly for
$\tu\in \tG^F$ unipotent and for $\tu'\in\tM^F$ unipotent we
set $\g_{\tM,\ti V}^{\tG}(\tu,\tu')=\sum_i(-1)^i\tr((\tu,\tu'),H^i_c(\tX,\bbq))$. Assuming that
$u=\t(\tu),u'=\t(\tu')$ we show:

(a) $\g_{\tM}^{\tG}(\tu,\tu')=|K^F|\g_M^G(u,v)$.
\nl
The restriction of $\t:\tG@>>>G$ defines a map $\tX@>>>X$ which is a principal covering with group $K^F$.
Hence we can identify $H^i_c(X,\bbq)$ with the $K^F$-invariants in $H^i_c(\tX,\bbq)$. It follows that (a) can be 
restated as follows:

$\sum_i(-1)^i\tr((\tu,\tu'),H^i_c(\tX,\bbq))=\sum_{k\in K^F}\tr((k\tu,\tu'),H^i_c(\tX,\bbq))$.
\nl
By the fixed point formula \cite{DL,3.2} we have $\tr((k\tu,\tu'),H^i_c(\tX,\bbq))=0$ for any $k\in K^F-\{1\}$
(since the fixed point of translation by $k$ on $\tX$ is empty). The desired equality follows.

We shall now omit the symbol $V$ in $\g_{M,V}^G(u,u')$; we write instead $\g_M^G(u,u')$.

\subhead 2.2\endsubhead
In the setup of 2.1 we define $a:cl(\tG^F)@>>>cl(G^F)$ and $a':cl(\tM^F)@>>>cl(M^F)$ by  
$$(a\ti f)(g)=\sum_{h\in\tG^F;\t(h)=g}\ti f(h),$$
$$(a'\ti f)(g)=\sum_{h\in\tM^F;\t(h)=g}\ti f(h).$$ 
For any $f\in cl(\tM^F)$ we show:
$$a(R_{\tM,\tP}^{\tG}(f))= R_{M,P}^G(a'(f)).\tag a$$
We must show that for $g\in G^F$ we have
$$\sum_{h\in\tG^F;\t(h)=g}R_{\tM}^{\tG}(f)(h)= R_M^G(a'(f))(g)$$
or (using \cite{L90, 1.7(b)}) that 
$$\align&\sum\Sb h\in\tG^F;\\\t(h)=g\endSb|\tM^F|\i|\tG_{h_s}^{0F}|\i
\sum\Sb z\in\tG^F;\\z\i h_sz\in\tM\endSb\sum\Sb \tv\in z\tM z\i\cap\tG_{h_s}^{0F}\\;unip.\endSb
\g_{z\tM z\i\cap\tG_{h_s}^0}^{\tG_{h_s}^0}(h_u,\tv)f(z\i h_s\tv z)\\&
=|M^F|\i|G_{g_s}^{0F}|\i\sum\Sb x\in G^F;\\x\i g_sx\in M\endSb\sum\Sb v\in xMx\i\cap G_{g_s}^{0F};\\unip.\endSb
\sum\Sb\tm\in\tM^F;\\\t(\tm)=x\i g_svx\endSb\g_{xMx\i\cap G_{g_s}^0}^{G_{g_s}^0}(g_u,v)f(\tm).\tag b\endalign$$
The right hand side of (b) is
$$\align&|K^F|\i\sum\Sb h\in\tG^F;\\\t(h)=g\endSb|\tM^F|\i|\tG_{h_s}^{0F}|\i|K_F|^2
|K^F|\i\sum\Sb z\in\tG^F;\\z\i h_sz\in\tM\endSb\sum\Sb\tv\in z\tM z\i\cap\tG_{h_s}^{0F};\\unip.\endSb
\\&
\sum\Sb\tm\in\tM^F;\\\t(\tm)=\t(z\i h_s\tv z)\endSb|K^F|\i\g_{z\tM z\i\cap\tG_{h_s}^0}^{\tG_{h_s}^0}(h_u,\tv)f(\tm)\\&
=|K^F|\i\sum\Sb h\in\tG^F;\\\t(h)=g\endSb|\tM^F|\i|\tG_{h_s}^{0F}|\i|K_F|^2
|K^F|\i\\&\sum\Sb z\in\tG^F;\\z\i h_sz\in\tM\endSb\sum\Sb\tv\in z\tM z\i\cap\tG_{h_s}^{0F};\\unip.\endSb
\sum_{k\in K^F}|K^F|\i\g_{z\tM z\i\cap\tG_{h_s}^0}^{\tG_{h_s}^0}(h_u,\tv)f(kz\i h_s\tv z).\endalign$$
(We have used 2.1(a).) This is the same as the left hand side of (b). This proves (a).

\subhead 2.3\endsubhead
We prove 1.9(b) for $G$ instead of $M$.
Let $\Th_G$ be the set of all pairs $(D,\cx)$ where $D$ is a conjugacy class of $G$ and $\cx\in ls_G(D)$
is irreducible (up to isomorphism). Now $F$ acts on $\Th_G$ by $F(D,\cx)=(FD,F^*\cx)$.
For $(D,\cx)\in \Th_G^F$ we denote by $\cl_{D,\cx}$ the line in $cl(G^F)$ containing the function
$\c_{\cx^\sha,\ph^\sha}:G^F@>>>\bbq$ where $\ph:F^*\cx@>\si>>\cx$ is an isomorphism; note that
$\c_{\cx^\sha,\ph^\sha}$ is equal to $0$ outside the closure of $D$. (This line is well defined.)
It is well known and easy to see that 

$cl(G^F)=\op_{(D,\cx)\in\Th_G^F}\cl_{D,\cx}$.
\nl
Hence to prove that $cl'(G^F)=cl(G^F)$ it is enough to show that 

(a) if $(D,\cx)\in\Th_G^F$ and $\ph:F^*\cx@>\si>>\cx$, $f_0=\c_{\cx^\sha,\ph^\sha}$, then $f_0\in cl'(G^F)$.
\nl
In the special case where $G_{der}$ is simply connected, this follows from \cite{L90, 9.5}. We shall 
deduce the general case from this special case. We can find $\t:\tG@>>>G,F:\tG@>>>\tG,K$ as in 2.1 such
that $\tG_{der}$ is simply connected.
Let $a:cl(\tG^F)@>>>cl(G^F)$ be as in 2.2. We define a linear map $b:cl(G^F)@>>>cl(\tG^F)$ by 
$(bf)(\tg)=f(\t(\tg))$; for $f\in cl(G^F)$ we have $abf=|K^F|f$.
Since 1.9(b) holds for $\tG$, we have $bf_0\in\cl'(\tG^F)$ hence $|K^F|f_0=abf_0\in a(\cl'(\tG^F))$. Thus it 
is enough to show that $a(\cl'(\tG^F))\sub cl'(G^F)$.

Let $(\tL,\tS)\in A_{\tG}$ be such that $F(\tL)=\tL,F(\tS)=\tS$ and let $\cf\in ls_{\tL}(\tS)$ be irreducible 
cuspidal with a given isomorphism $\ps:F^*\cf@>\si>>\cf$. It is enough to show that

(b) $a(R_{\tL}^{\tG}(\c_{\cf^\sha,\ps^\sha}))\in cl'(G^F)$.
\nl
Let $L=\t(\tL),S=\t(\tS)$; we have $(L,S)\in A_M$. Let $\t':\tL@>>>L$ be the restriction of $\t$; we define
define $a':cl(\tL^F)@>>>cl(L^F)$ by $(a'\ti f)(g)=\sum_{\tg\in\tL^F;\t'(\tg)=g}\ti f(\tg)$. By 2.2(a), for any 
$f\in cl(\tL^F)$ we have

(c) $a(R_{\tL}^{\tG}(f))= R_L^G(a'(f))$.
\nl
From this we see that the left hand side of (b) is equal to $R_L^G(a'(\c_{\cf^\sha,\ps^\sha}))$. From the 
definitions we see that $a'(\c_{\cf^\sha,\ps^\sha})$ is a linear combination of functions of the form
$\c_{\ce_0^\sha,\ph_0^\sha}:L^F@>>>\bbq$ where $\ce_0\in ls_L(S)$ is irreducible cuspidal and 
$\ph_0:F^*\ce_0@>\si>>\ce_0$. It follows that $R_L^G(a'(\c_{\cf^\sha,\ps^\sha}))\in cl'(G^F)$. We 
see that (b) holds. This completes the proof of 1.9(b) for $G$.

\head 3. A direct sum decomposition of $cl(G^F)$\endhead
\subhead 3.1\endsubhead
In this section there is no restriction on $q$.
Let $(D,\cx)\in\Th_G$. We associate to $(D,\cx)$ an admissible stratum of $G$. Let $E$ be the set of semisimple 
parts of elements in $D$; this is a conjugacy class in $G$. For $s\in E$ let $[s]$ be the set of unipotent 
conjugacy classes of $G_s^0$ such that $sC\sub D$. For any $s\in E$ and $C\in[s]$ we define $f_s:C@>>>D$ by 
$u\m su$; then $f_s^*\cx\in ls_{G_s^0}(C)$. Let $fs^*\cx=\op_{\cy\in Q_{s,C}}\cy$ be
the isotypic decomposition of $f_s^*\cx$; thus each $\cy$ is an isotypic object of $ls_{G_s^0}(C)$. 
Let $D'$ be the set of all pairs $(g,\cy)$ where $g\in D$ and $\cy\in Q_{g_s,C}$ where $C\in[g_s]$ contains
$g_u$. Then $D'$ is naturally an algebraic variety with a transitive action of $G$ such that the 
map $D'@>>>D$, $(g,\cy)\m g$ is a $G$-equivariant unramified finite covering. For $s\in E$, $C\in[s]$, 
$\cy\in Q_{s,C}$, we choose an irreducible summand $\et$ of $\cy$; the generalized Springer correspondence 
\cite{L84, 6.3} for the reductive connected group $G_s^0$ associates to the pair 
$(C,\et)$ a triple $(L,S,\cf)=(L_\cy,S_\cy,\cf_\cy)$ (up to $G_s^0$-conjugacy) where $L\in\L(G_s^0)$, 
$S=\cz_L^0c$ with $c=c_\cy$ a unipotent class of $L$ and $\cf=\bbq\bxt\cf_0\in ls_L(S)$ is irreducible cuspidal 
with $\cf_0\in ls_L(c)$ irreducible; this triple is independent of the choice of $\et$ since $\cy$ is isotypic. 
Let $M=M_\cy=Z_G(\cz_L^0)\in\L(G)$. Let $D_\cy$ be the conjugacy class in $M$ containing $sc$. Let 
$\Si=\Si_\cy=D_\cy\cz_M^0$. Since $L\in\L(G_s^0)$, we have $Z_{G_s^0}(\cz^0_L)=L$ hence $(Z_{G_s}(\cz^0_L))^0=L$.
We have $M_s=G_s\cap M=G_s\cap Z_G(\cz_L^0)=Z_{G_s}(\cz^0_L)$ so that $M_s^0=L$. We have 
$Z_M(\cz^0_{M_s^0})=Z_M(\cz^0_L)=Z_G(\cz^0_L)\cap M=M$ hence $s$ is isolated in $M$ and $\Si$ is an isolated 
stratum of $M$. Hence we can define $Y=Y^G_{M,\Si}$, a stratum of $G$. If $(L,S,\cf)$ is replaced by a 
$G_s^0$-conjugate or if $(s,C,\cy)$ is replaced by a triple in the same $G$-orbit, then $Y$ is replaced by a 
$G$-conjugate hence it remains the same. Thus the stratum $Y$ depends only on $(D,\cx)$.
For $\cy,(L,S,\cf),M,\Si$ as above we can find $\cf'\in ls_M(\Si)$ irreducible such that
the inverse image of $\cf'$ under $C@>>>\Si,u\m su$ contains $\cf$ as a direct summand. By the arguments
in \cite{L84, 2.10}, $\cf'$ is cuspidal. It follows that $Y$ is an admissible stratum. We set $Y=\ps(D,\cx)$.

Note that if $(D,\cx)\in\Th_G^F$ then $F(Y)=Y$.

\subhead 3.2\endsubhead
Let $\G'_G$ be the set of all triples $(L,S,\ce_0)$ where $(L,S)\in A_G$ is such that $FL=L,FS=S$ and
$\ce_0\in ls_L(S)$ is irreducible cuspidal (up to isomorphism) such that
$F^*\ce_0\cong\ce_0$. Let $\G_G$ be the set of orbits of the conjugation action of $G^F$ in $\G'_G$.
For $(L,S,\ce_0)\in\G'_G$ we choose an  isomorphism $\ph_0:F^*\ce_0@>\si>>\ce_0$ of local systems on $S$. Then
$\c_{\ce_0^\sha,\ph_0^\sha}:L^F@>>>\bbq$ is well defined; it is a class function on $L^F$.
Let $\cl_{L,S,\ce_0}$ be the line in $cl(G^F)$ containing $\car_L^G(\c_{\ce_0^\sha,\ph_0^\sha})$ for some/any 
$\ph_0$ as above; this line depends only the image of $(L,S,\ce_0)$ in $\G_G$. We have the following result.

\proclaim{Theorem 3.3} (i) We have $cl(G^F)=\op_{(L,S,\ce_0)\in\G_G}\cl_{L,S,\ce_0}$.

(ii) For any $F$-stable admissible stratum $Y$ of $G$ we define $cl_Y(G^F)$ to be the subspace
$\sum_{(L,S,\ce_0)\in\G_G;Y^G_{L,S}=Y}\cl_{L,S,\ce_0}$ of $cl(G^F)$ (this is a direct sum, see (i)); 
we define $\un{cl}_Y(G^F)$ to be 
the subspace $\op_{(D,\cx)\in\Th_G^F;\ps(D,\cx)=Y}\cl_{D,\cx}$ of $cl(G^F)$ (see 2.3, 3.1). We have
$cl_Y(G^F)=\un{cl}_Y(G^F)$ and $cl(G^F)=\op_Ycl_Y(G^F)$ where 
$Y$ runs over the $F$-stable admissible strata of $G$. 
\endproclaim
The fact that the sum in (i) is direct follows from the orthogonality relations \cite{L85, 9.9} (its hypotheses are
satisfied by the results in \cite{L86} and \cite{L12}). If $(D,\cx)\in \Th_G^F$ and $Y=\ps(D,\cx)$ then we have

(a) $\cl_{D,\cx}\sub cl_Y(G^F)$.
\nl
When $G_{der}$ is simply connected, (a) follows from \cite{L90, 9.5}. (One can replace $R_L^G$ in 
{\it loc.cit.} with $q$ large by $\car_L^G$ without restriction on $q$.) The general case can be reduced to this
special case by passage to $\tG$ as in the proof in 2.3 (again replacing $R_L^G$ by $\car_L^G$).
Since the lines $\cl_{D,\cx}$ span $cl(G^F)$ we see that (a) implies that the sum in (i) is equal to $cl(G^F)$.
Thus (i) holds. From (a) we see that $\un{cl}_Y(G^F)\sub cl_Y(G^F)$ for any $Y$. Since
$\op_Y\un{cl}_Y(G^F)=\op_Y cl_Y(G^F)=cl(G^F)$ (see 2.3 and (i)) it follows that $\un{cl}_Y(G^F)=cl_Y(G^F)$ for 
any $Y$. This proves (ii).

\subhead 3.4\endsubhead
From 3.3 and the orthogonality relations mentioned in the proof of 3.3 one can deduce that the ``Mackey formula''
for $R_{L,P}^G$ stated by Deligne (unpublished) in 1976 for $q$ large and in \cite{BM} for $q>2$ remains valid 
without restriction on $q$ if $R_{L,P}^G$ is replaced by $\car_L^G$.

\subhead 3.5\endsubhead
Let $(D,\cx)\in\Th_G$. We use notation of 3.1. We say that $(D,\cx)$ is of principal type if for $s\in E$,
$C\in[s]$, the local system $f_s^*\cx$ on $C$ is such that some/any irreducible summand $\eta$ of $f_s^*\cx$ is
such that $(C,\eta)$ appears in the usual Springer correspondence for $G_s^0$. An equivalent condition is that 
the stratum $Y=\ps(D,\cx)$ is the variety of regular semisimple elements in $G$.
For example, $(D,\bbq)$ is of principal type.

Now let $(D,\cx)\in\Th_G^F$ be of principal type; let $\ph:F^*\cx@>\si>>\cx$ be an isomorphism. From 3.3(a) we 
deduce
$$\c_{\cx^\sha,\ph^\sha}=\sum_{T,\th}c_{D,\cx;T,\th}\car_T^G(\th)\tag a$$
where $T$ runs over the $F$-stable maximal tori in $G$, $\th$ runs through the set
of characters $T^F@>>>\bbq$ and the pairs $(T,\th)$ are taken up to $G^F$-conjugacy; $c_{D,\cx;T,\th}\in\bbq$ are 
uniquely determined. Equivalently, we have
$$\c_{\cx^\sha,\ph^\sha}=\sum_{T,\th}c_{D,\cx;T,\th}R_T^G(\th).\tag b$$
Indeed, we have $R_T^G(\th(=\car_T^G(\th)$. This follows from the results in \cite{L90} (for large $q$) and their
extension to general $q$ in \cite{Sh}. 
Moreover, from \cite{L90, 9.5} we see that $c_{D,\cx;T,\th}$ are explicitly known (at least 
if $G_{der}$ is simply connected, but the general case can be reduced to this case as before).
Since the multiplicities of irreducible representations of $G^F$ in $R_T^G(\th)$ are known, it follows that the 
functions $\c_{\cx^\sha,\ph^\sha}$ are computable as explicit linear combinations of irreducible characters.

In particular, (a),(b) hold when $D$ is an $F$-stable conjugacy class in $G^F$ and $\cx=\bbq$.

\subhead 3.6\endsubhead
Let $(D,\cx)\in\Th^F_G$. Let $\cz$ be the set of all $(D',\cx')\in\Th^F_G$ such that $D'$ is contained in the 
closure of $D$. For any $(D',\cx')\in\cz$ we choose an isomorphism $\ph_{\cx'}:F^*\cx'@>\si>>\cx'$. We have 

(a) $\c_{\cx,\ph_{\cx}}=\sum_{(D',\cx')\in\cz}d_{D',\cx'}\c_{\cx'{}^\sha,\ph_{\cx'}^\sha}$
\nl
where $d_{D',\cx'}\in\bbq$. Assume now that $(D,\cx)$ is of principal type. Then $d_{D',\cx'}=0$ unless
$(D',\cx')$ is of principal type. (This can be deduced from the results in \cite{L86} on Green functions.)
Using 3.5(b) we deduce
$$\c_{\cx,\ph}=\sum_{T,\th}\ti c_{D,\cx;T,\th}R_T^G(\th)\tag b$$
where $\ti c_{D,\cx;T,\th}\in\bbq$ is explicitly computable.
In particular, (b) holds when $\cx=\bbq$. We see that:

(c) {\it the class function on $G^F$ equal to $1$ on $D^F$ and equal to $0$ on $G^F-D^F$ is a linear combination
of functions of the form $R_T^\th$.}
\nl
This has been conjectured in \cite{L78, 2.16}. Note that the coefficients in the linear combination above are 
explicitly computable. Since each $R_T^\th$ is an explicit linear combination of irreducible characters, we 
deduce that for any $D$ as above the average value on $D^F$ of any irreducible character of $G^F$ is explicitly 
computable. In the case where $D$ is a semisimple class, a result like (c) appears (in a stronger form) in
\cite{DL, 7.5}.

\enddocument